\documentclass[12pt]{article}\usepackage{amsfonts,amssymb}

\usepackage{amsmath,bbm}
\usepackage{theorem}          
\usepackage{latexsym}
\usepackage{float}                   
\usepackage{times}
\usepackage{cite}

\newtheorem{tht}{Theorem}[section]

\newtheorem{thl}[tht]{Lemma}

\newcommand{\ang}{\raisebox{0.2ex}{\scriptsize$\triangleright$}}

\newcommand{\mn}{\medskip}    

\newcommand{\rti}{\,{\scriptstyle\rtimes}\,} %

\newcommand{\cG}{{\mathcal{G}}} 
\newcommand{\Hh}{{\mathcal{H}}} 
\newcommand{\cO}{{\mathcal{O}}}
\newcommand{\cK}{{\mathcal{K}}}
\newcommand{\cU}{{\mathcal{U}}}   

\newcommand{\cX}{{\mathcal{X}}}

\newcommand{\cS}{\mathcal{S}}

\newcommand{\cY}{\mathcal{Y}}

\newcommand{\dR}{\mathbb{R}}

\newcommand{\dN}{\mathbb{N}}

\newcommand{\Lin}{{\mathrm{Lin}}}

\newcommand{\Ker}{\mathrm{ker}\,}

\newcommand{\SU}{\cO(\mathrm{SU}_q(2))}

\newcommand{\su}{\cU_q(\mathrm{su}_2)}

\begin{document}

\date{\small{Fakult\"at f\"ur Mathematik und 
Informatik\\ Universit\"at Leipzig, 
Augustusplatz 10, 04109 Leipzig, Germany\\ 
E-mail: schmuedg@mathematik.uni-leipzig.de / 
wagner@mathematik.uni-leipzig.de}
}

\title{Operator representations of cross product algebras 
of Podles' quantum spheres}

\author{Konrad Schm\"udgen  and Elmar Wagner}

\maketitle

\renewcommand{\theenumi}{\roman{enumi}}
\begin{abstract}
Operator representations of the 
cross product $\ast$-al\-ge\-bra
 $\cO(S^2_{qc})\rti\cU_q(\mathrm{su}_2)$ of the Hopf $\ast$-al\-ge\-bra 
$\cU_q(\mathrm{su}_2)$ and 
its module $\ast$-al\-ge\-bras $\cO(S^2_{qc})$ of Podles' spheres are studied. 
Two classes of representations are described by explicit formulas for the 
actions of the generators.
\end{abstract}
%
%
%
\setcounter{section}{-1}
\section{Introduction}
In this paper we investigate Hilbert space representations of the left 
cross pro\-duct $\ast$-al\-ge\-bras $\cO(S^2_{qc})\rti \cU_q(\mathrm{su}_2)$ 
of the 
quantized enveloping algebras $\cU_q(\mathrm{su}_2)$ 
and the coordinate algebras 
$\cO(S^2_{qc})$, $c\in [0,+\infty]$, of  the Podles' spheres. As discussed 
in the introduction of our previous paper \cite{SW}, there are two principal 
methods. In the first approach, we begin with a representation of the 
$\ast$-al\-ge\-bras $\cO(S^2_{qc})$ given in a canonical form and extend 
it to a representation of $\cO(S^2_{qc})\rti\cU_q(\mathrm{su}_2)$. 
The main technical 
tool for this is to "decouple'' the cross relations of the cross product 
algebra by finding an auxiliary $\ast$-sub\-al\-ge\-bra $\cY_c$ of a slightly 
larger cross product $\ast$-al\-ge\-bra 
$\hat{\cO} (S^2_{qc})\rti\cU_q(\mathrm{su}_2)$ 
which commutes with the $\ast$-al\-ge\-bra $\cO(S^2_{qc})$. Such decouplings 
have been found and studied in \cite{F}. Since elements of $\cY_c$ act 
as unbounded operators, we impose some regularity conditions
during the derivation of the corresponding representations. 
In the second approach, we assume 
that there 
exists some 
$l \in \frac{1}{2}\dN_0$ such 
that the restriction  to $\cU_q(\mathrm{su}_2)$ is 
a direct sum of representations $T_{l+n}$, $n\in\dN_0$. 
The Heisenberg representation (see \cite{SW}, 5.2) of 
$\cO(S^2_{qc})\rti\cU_q(\mathrm{su}_2)$ satisfies this condition with $l=0$. 
We show that for any $l\in \frac{1}{2}\dN$ there exist precisely 
two inequivalent irreducible representations of 
$\cO(S^2_{qc})\rti \cU_q(\mathrm{su}_2)$ having this property. 
In the case $l=1/2$ the existence of these representations was known 
for two special quantum spheres. 
For the Podles' sphere $S^2_{q,\infty}$ related to the quantum vector 
space $\dR^3_q$ they appeared in \cite{SW}, Section 6.5. 
For the standard quantum sphere $S^2_{q,0}$ the
representations with $l=1/2$ are contained in \cite{DS}.

This paper is organized as follows. The definition of the cross product 
algebra $\cO(S^2_{qc})\rti\cU_q(\mathrm{su}_2)$ is given in Section 1. 
The "decoupling'' procedure is developed in Section 2. Representations 
of the auxiliary $\ast$-al\-ge\-bra $\cY_c$ are studied in Section 3. 
The first approach as explained above is carried out in Section 4. 
The second approach is contained in Section 5.

All facts and notions on quantum groups used in this paper can be found, 
for instance, in \cite{KS}. The algebra $\cU_q(\mathrm{sl}_2)$ and 
the Podles' spheres $S^2_{qc}$ have been discovered in \cite{KR} 
and \cite{P}, respectively.

Throughout this paper, $q$ stands for a real number of the open interval 
$(0,1)$. We abbreviate $\lambda:=(q{-}q^{-1})$, 
$\lambda_n:=(1{-}q^{2n})^{1/2}$ and 
$[n]:=(q^n{-}q^{-n})/(q{-}q^{-1})$. 
%
\section{Definition of the cross product algebra}
Let $\cU$ be a Hopf $\ast$-al\-ge\-bra and let $\cX$ be a left $\cU$-module 
$\ast$-al\-ge\-bra, that is,  $\cX$  is a unital $\ast$-al\-ge\-bra with left 
$\cU$-action $\ang$ satisfying 
$f\ang xy=(f_{(1)}\ang x) (f_{(2)}\ang y)$, $f\ang 1  =\varepsilon(f) 1$ 
and $(f\ang x)^\ast = S(f)^\ast \ang x^\ast $ for $x,y\in\cX$ and $f\in\cU$. 
Here $\Delta (f)=f_{(1)} \otimes f_{(2)}$ is the Sweedler notation for 
the comultiplication $\Delta(f)$ of $f\in \cU$. 
Then the {\it left cross product $\ast$-al\-ge\-bra} $\cX\rti\cU$ is 
the $\ast$-al\-ge\-bra generated by the two $\ast$-sub\-al\-ge\-bras $\cX$ 
and $\cU$ with respect to the cross commutation relation 
\begin{equation}\label{fxrel}
fx = (f_{(1)}\ang x) f_{(2)}, \ \ x\in\cX,\ f\in\cU.
\end{equation}
The Hopf $\ast$-al\-ge\-bra $\cU_q(\mathrm{su}_2)$ is generated by elements 
$E,F,K,K^{-1}$ with relations
\begin{equation}\label{urel}
KK^{-1} \!=\! K^{-1} K\!=\!1,\, KE\!=\!qEK,\, FK\!=\!qKF,\,
EF- FE\!=\!\lambda^{-1}(K^2-K^{-2}),
\end{equation}
involution $E^\ast=F$, $K=K$ and comultiplication 
$$
\Delta (E)\!=\!E\otimes K+K^{-1} \otimes E,\
\Delta (F)\!=\!F\otimes K + K^{-1}\otimes K,\  \Delta (K)\!=\!K\otimes K.
$$
There is a dual pairing $\langle \cdot,\cdot\rangle$ of Hopf 
$\ast$-al\-ge\-bras $\su$ and $\SU$ given on generators by
$$
\langle K^{\pm 1},d\rangle = \langle K^{\mp 1},a\rangle = q^{\pm 1/2}, \quad
\langle E,c\rangle = \langle F,b\rangle=1
$$
and zero otherwise, where $a,b,c,d$ are the usual generators 
of the Hopf algebra $\SU$ (see e.g. \cite{KS}, Chapter 4). 

We shall use the definition of the coordinate algebras 
$\cO(S^2_{qc})$, $c\in [0,+\infty]$, of Podles' spheres as given in \cite{P}. 
For $c\in[0,+\infty)$, $\cO(S^2_{qc})$ is the $\ast$-al\-ge\-bra 
with generators $A\!=\!A^\ast, B,B^\ast$ and defining relations
\begin{equation}\label{podrel}
AB\!=\!q^{-2} BA,\ AB^\ast\!=\!q^2 B^\ast A,\ B^\ast B\!=\!A-A^2+c,\ 
BB^\ast \!=\!q^2 A-q^4 A^2+c.
\end{equation}
For $c=+\infty$, the defining relations of $\cO(S^2_{q ,+\infty})$ are
\begin{equation}\label{prodrel1}
AB=q^{-2} BA,\  AB^\ast = q^2B^\ast A,\ B^\ast B=-A^2+1,\ BB^\ast=-q^4 A^2+1.
\end{equation}
Let $c<\infty$. 
As shown in \cite{P}, 
$\cO(S^2_{qc})$ is a right $\SU$-comodule $\ast$-al\-ge\-bra such that 
\begin{equation}\label{xgener}
x_{-1}:= q^{-1}(1+q^{2})^{1/2}B,\ \, x_1 :=-(1+q^{2})^{1/2}B^\ast,\ \,
x_0 := 1-(1+q^2)A.
\end{equation}
transform by the spin $1$ matrix corepresentation 
$(t^1_{ij})$ of $\mathrm{SU}_q(2)$. Hence $\cO(S^2_{qc})$ 
is a left $\su$-module 
$\ast$-al\-ge\-bra with left action given by 
$f\ang x_j=\sum_i x_i\langle f,t^1_{ij}\rangle$ 
for $f\in \su$, $j=-1,0,1$. Inserting the form of the 
matrix $(t^1_{ij})$ (see \cite{P}, p.\ 194, or \cite{KS}, p.\ 124) and 
the Hopf algebra pairing $\langle\cdot,\cdot\rangle$ into (\ref{fxrel})  
we derive the following cross relations for the cross product algebra
$\cO(S^2_{qc})\rti \su$:
\begin{align*}
&KA \!=\! AK,\  EA \!=\! AE + q^{-1/2} B^\ast K,\ 
FA \!=\!AF - q^{-3/2} BK,\\
&KB \!=\! q^{-1} BK,\ EB \!=\! q BE -q^{1/2} (1+q^2) AK + q^{1/2} K, \  
FB \!=\! q BF,\\
&KB^\ast \!=\! q B^\ast K,\, EB^\ast \!=\! q^{-1} B^\ast E, \,
FB^\ast \!=\! q^{-1} B^\ast F + q^{-1/2} (1+q^2) AK - q^{-1/2} K.
\end{align*}
That is, $\cO(S^2_{qc})\rti \cU_q(\mathrm{su}_2)$ is the algebra 
with generators $A$, $B$, $B^\ast$, $E$, $F$, $K$, $K^{-1}$, with defining 
relations (\ref{urel}), (\ref{podrel}) (resp.\ (\ref{prodrel1})) and 
the preceding set of cross relations.

For $c=\infty$, we set 
\begin{equation}\label{xgener1}
x_{-1}:= q^{-1}(1+q^{2})^{1/2}B,\ \,x_1 :=-(1+q^{2})^{1/2}B^\ast,\ \,
x_0 := -(1+q^2)A.
\end{equation}
Then the cross relations for the cross product algebra
$\cO(S^2_{q\infty})\rti \su$ can be written as 
\begin{align*}
&KA = AK,\ \, EA = AE + q^{-1/2} B^\ast K,\ \,
FA =AF - q^{-3/2} BK,\\
&KB = q^{-1} BK,\ \, EB = q BE -q^{1/2} (1+q^2) AK, \ \,  
FB = q BF,\\
&KB^\ast = q B^\ast K,\,\ EB^\ast = q^{-1} B^\ast E, \, \ 
FB^\ast = q^{-1} B^\ast F + q^{-1/2} (1+q^2) AK.
\end{align*}

\section{"Decoupling'' of the cross product algebra}
Let us first suppose that $c\in[0,+\infty)$.
From the relations $AB=q^{-2} BA$ and $AB^\ast = q^2 B^\ast A$ it is 
clear that $\cS=\{A^n\,;\, n\in\dN_0\}$ is a left and right Ore set 
of the algebra $\cO(S^2_{qc})$. Moreover the algebra $\cO(S^2_{qc})$ has 
no zero divisors. Hence the localization algebra, 
denoted by $\hat{\cO}(S^2_{qc})$,  of $\cO(S^2_{qc})$ at $\cS$ exists. 
The $\ast$-al\-ge\-bra $\cO(S^2_{qc})$ is then a $\ast$-sub\-al\-ge\-bra 
of $\hat{\cO}(S^2_{qc})$ and all elements $A^n$, $n\!\in\!\dN_0$, 
are invertible in $\hat{\cO}(S^2_{qc})$. From Theorem 3.4.1 in \cite{LR} 
we conclude that $\hat{\cO} (S^2_{qc})$ is a left (resp.\ right) 
$\su$-module $\ast$-al\-ge\-bra which contains $\cO(S^2_{qc})$ as a left 
(resp.\ right) $\su$-module $\ast$-sub\-al\-ge\-bra. The left action of the 
generators $E,F,K$ on $A^{-1}$ is given by
$$
E\ang A^{-1} = -q^{-5/2} B^\ast A^{-2},\ \, F\ang A^{-1} = q^{1/2} BA^{-2}, 
\ \,K\ang A^{-1}=A^{-1}.
$$
Hence the left cross product algebra $\hat{\cO}(S^2_{qc})\rti \su$ is a 
well-defined $\ast$-al\-ge\-bra.

Let $\cY_c$ denote the $\ast$-sub\-al\-ge\-bra of $\hat{\cO}(S^2_{qc})\rti\su$ 
generated by
\begin{align}\label{xydef}
&X := q^{3/2} \lambda FK^{-1} A + q B, \ \,
X^\ast = q^{3/2} \lambda AK^{-1} E + q B^\ast, \\
&Y := qK^{-2} A,\ \,Y^{-1} = q^{-1} A^{-1}  K^2.\label{xydef1}
\end{align}
Note that $Y=Y^\ast$. It is straightforward to check that 
the elements $X$, $X^\ast$, $Y$, $Y^{-1}$ commute with the generators 
$A$, $B$, $B^\ast$ of $\cO(S^2_{qc})$. Hence the algebras $\cY_c$ and 
$\hat{\cO}(S^2_{qc})$ {\it commute} inside the cross product algebra 
$\hat{\cO} (S^2_{qc})\rti \su$. Moreover, the generators of 
$\cY_c$ satisfy the commutation relations 
\begin{equation}\label{xyrel}
YX = q^2 XY,\ \, YX^\ast = q^{-2} X^\ast Y,\ \,
X^\ast X {-} q^2 XX^\ast = (1-q^2) (Y^2+c).
\end{equation}

We denote by $\cU^\prime_q(\mathrm{su}_2)$ the Hopf $\ast$-sub\-al\-ge\-bra 
generated by 
$e:= EK$, $f:= K^{-1} F$  and $k:=K^2.$ As an algebra, 
$\cU^\prime_q(\mathrm{su}_2)$ has generators $e$, $f$, $k$, $k^{-1}$ 
with defining 
relations 
\begin{equation}\label{efrel}
kk^{-1} = k^{-1}k=1,\ \, ke = q^2 ek,\ \,kf = q^{-2} fk,\ \,
ef {-} fe = \lambda^{-1} (k{-}k^{-1}).
\end{equation}
The comultiplication on the generators $e$, $f$, $k$ is given by 
$$
\Delta(e)= e\otimes k + 1\otimes e, 
\ \,\Delta (f) = f\otimes 1 + k^{-1} \otimes f,\ \, \Delta (k)=k\otimes k.
$$
From (\ref{xydef}) and (\ref{xydef1}) we obtain
\begin{equation}\label{efdef}
f \!=\! q^{-1/2} \lambda^{-1} (X{-}qB) A^{-1},\ \, 
e \!=\! q^{1/2} \lambda^{-1} (X^\ast{-} q^{-1} B^\ast) Y^{-1},\ \, 
k\!=\!qY^{-1} A.
\end{equation}
Hence the two commuting algebras $\cY_c$ and $\hat{\cO}(S^2_{qc})$ 
generate the $\ast$-sub\-al\-ge\-bra 
$\hat{\cO}(S^2_{qc})\rti \cU^\prime_q (\mathrm{su}_2)$ 
of $\hat{\cO}(S^2_{qc})\rti \cU_q(\mathrm{su}_2)$. 
So one can say that $\cY_c$ and $\hat{\cO}(S^2_{qc})$ 
completely "decouple'' the cross product algebra 
$\hat{\cO}(S^2_{qc})\rti \cU^\prime_q (\mathrm{su}_2)$. 

There is an alternative way to define the $\ast$-al\-ge\-bra 
$\hat{\cO}(S^2_{qc})\rti \cU^\prime_q (\mathrm{su}_2)$ 
by taking the two sets $A\!=\!A^\ast$, $B$, $B^\ast$ and 
$X$, $X^\ast$, $Y\!=\!Y^\ast$, $Y^{-1}$ of 
pairwise commuting generators 
with defining relations (\ref{podrel}), (\ref{xyrel}) 
and the obvious relations 
\begin{equation}\label{obrel}
AA^{-1} = A^{-1} A=1, \quad  YY^{-1} = Y^{-1} Y=1.
\end{equation}
Indeed, if we define $e,f,k$ by (\ref{efdef}), then the 
relations (\ref{efrel}) of $\cU_q(\mathrm{su}_2)$ and the cross relations 
of $\hat{\cO}(S^2_{qc})\rti \cU^\prime_q(\mathrm{su}_2)$ can be derived 
from this set of defining relations. 

The larger cross product algebra 
$\hat{\cO}(S^2_{qc})\rti \cU_q(\mathrm{su}_2)$ can be redefined in 
a similar manner 
if we replace the generator $Y$ by $K$. That is, 
$\hat{\cO}(S^2_{qc})\rti \cU_q(\mathrm{su}_2)$ is the $\ast$-al\-ge\-bra 
with the two sets $A\!=\!A^\ast$, $B$, $B^\ast$ 
and $X$, $X^\ast$, $K$, $K^{-1}$ of 
pairwise commuting generators 
with defining relations (\ref{podrel}), (\ref{obrel}), and 
\begin{equation}\label{xkrel}
KA\!=\!AK,\, BK\!=\!qKB ,\, KB^\ast \!=\! qB^\ast K,\, 
X^\ast X\!-\!q^2 XX^\ast \!=\! (1\!-\!q^2)(q^2 K^{-4} A^2\!+\!c).
\end{equation}
The generators $F$ and $E$ are then given by 
\begin{equation}\label{ferel}
F=q^{-3/2} \lambda^{-1} (X-qB) KA^{-1},\ \ 
E=q^{-3/2} \lambda^{-1} A^{-1} K(X^\ast-q B^\ast).
\end{equation}

The preceding considerations and facts carry over almost verbatim to 
the case $c=+\infty$. The only difference is that in the case $c=+\infty$ 
one has to set $c=1$ in   
the third equations of (\ref{xyrel}) and (\ref{xkrel}). 
%
\section{Operator representations of the $\ast$-al\-ge\-bra ${\bf \cY_c}$}
%
For the study of representations of the $\ast$-al\-ge\-bra $\cY_c$ we need two 
auxiliary lemmas. The first one restates the {\it Would decomposition} of 
an  isometry (see \cite{SF}, Theorem 1.1), while the second 
is Lemma 4.2 (ii) in \cite{SW}.

\begin{thl}                                                     \label{L1}
Each isometry $v$ on a Hilbert space $\cK$ is up to unitary equivalence 
of the following form: There exist Hilbert subspaces $\cK^u$ and $\cK^s_0$ 
of $\cK$ and a unitary operator $v_u$ on $\cK^u$ such that $v=v_u\oplus v_s$ 
on $\cK = \cK^u \oplus \cK^s$ and $v_s$ acts on 
$\cK^s = \oplus^\infty_{n=0} \cK^s_n$ by 
$v_s \zeta_n = \zeta_{n+1}$, where each $\cK^s_n$ is $\cK^s_0$ and 
$\zeta_n\in\cK^s_n$. 
Moreover, $\cK^u = \cap^\infty_{n=0} v^n \cK$.
\end{thl}

\begin{thl}                                                     \label{L2}
Let $v$ be the operator $v_s$ on $\cK^s = \oplus^\infty_{n=0} \cK^s_n$, 
$\cK^s_n= \cK^s_0$, from  Lemma \ref{L1} 
and let $Y$ be a self-adjoint operator 
on $\cK^s$ such that $q^2 v Y \subseteq Y v$. Then there is a self-adjoint 
operator $Y_0$ on the Hilbert space  $\cK^s_0$ such that 
$Y\zeta_n=q^{2n} Y_0 \zeta_n$ and $v\zeta_n=\zeta_{n+1}$ for 
$\zeta_n\in \cK^s_n$ and $n\in\dN_0$.
\end{thl}

We treat the case $c\in [0,+\infty)$.
Suppose that we have a representation of relations (\ref{xyrel}) by 
closed operators $X$, $X^\ast$ and a self-adjoint operator $Y$ with trivial 
kernel acting on a Hilbert space $\cK$.

Let $X= v|X|$ be the polar decomposition of the operator $X$. 
Since $q\in (0,1)$, $c\ge 0$ and $\ker Y=\{0\}$, we have $\ker X= \{0\}$ by 
the third equation of (\ref{xyrel}). Hence $v$ is an isometry on $\cK$, 
that is, $v^\ast v=I$.

By (\ref{xyrel}), we have $YX^\ast X= X^\ast XY$. We assume that the 
self-adjoint operators $Y$ and $X^\ast X$ strongly commute. Then $Y$ 
and $|X|=(X^\ast X)^{1/2}$ also strongly commute. 
Again by (\ref{xyrel}), 
$Y v|X|=YX = q^2 XY=q^2 v|X|Y =q^2 v Y|X|$. Since $\ker|X|=\ker X=\{0\}$, 
we conclude that
\begin{equation}\label{vyrel}
q^2 v Y=Yv.
\end{equation}
Since $X^\ast = |X| v^\ast$, the third equation of (\ref{xyrel}) rewrites as 
\begin{equation}\label{xvrel}
|X|^2 = q^2 v|X|^2 v^\ast + (1{-}q^2)(Y^2 +c).
\end{equation}
Multiplying by $v$ and using the relations $v^\ast v=I$ and (\ref{vyrel}) 
we get
$$
|X|^2 v= v\big(q^2 |X|^2 + (1{-}q^2)(q^4 Y^2 + c)\big).
$$
Proceeding by induction and using once more (\ref{vyrel}) we derive
\begin{equation}\label{xvnrel}
|X|^2 v^n = v^n \big(q^{2n} |X|^2 + (1{-}q^{2n})(q^{2n+2} Y^2 + c)\big).
\end{equation}
Now we use the Would decomposition $v=v_u\oplus v_s$ on 
$\cK=\cK^u\oplus \cK^s$ of the isometry $v$ by Lemma \ref{L1}. 
Since $\cK^u=\cap^\infty_{n=0} v^n \cK$ and $Yv=q^2 vY$, 
the Hilbert subspace $\cK^u$ reduces the self-adjoint operator $Y$. 
From (\ref{xvnrel}) it follows that $|X|^2$ leaves a dense subspace 
of the space $\cK^u$ invariant. We assume that $\cK^u$ reduces $|X|^2$. 
Hence we can consider all operators occurring in  relation (\ref{xvnrel}) 
on the subspace $\cK^u$. For the unitary part $v_u$ of $v$ we have 
$v_uv^\ast_u=I$ on $\cK^u$. Multiplying (\ref{xvnrel}) by 
$(v^\ast_u)^n$ and using again relation (\ref{vyrel}) we derive on $\cK^u$
$$
0\le q^{2n} v^n_u |X|^2 (v^\ast_u)^n 
= |X|^2 - (1{-}q^{2n})c + q^2 Y^2 - q^{-2n+2} Y^2
$$
for all $n\in \dN$. Letting $n\rightarrow\infty$ and remembering that  
$\ker Y=\{0\}$ and $q\in (0,1)$ we conclude that the latter is only possible 
when $\cK^u = \{0\}$. That is, the isometry $v$ 
is (unitarily equivalent to ) the unilateral shift operator $v_s$.

Since $v=v_s$, Lemma \ref{L2} applies to the relation $q^2 v Y=Yv$. 
Hence there exists a self-adjoint operator $Y_0$ on $\cK^s_0$ 
such that $Y\zeta_n = q^{2n} Y_0 \zeta_n$ 
on $\cK^s=\oplus^\infty_{n=0} \cK^s_n$, $\cK^s_n = \cK^s_0$. 
Since $v^\ast\zeta_0=0$, (\ref{xvrel}) yields 
$|X|^2\zeta_0= (1-q^2)(Y^2_0 +c)\zeta_0$. 
Using this equation and (\ref{xvrel}) we compute 
\begin{align}\label{x2rel}
|X|^2\zeta_n&=|X|^2 v^n\zeta_0 
= v^n \big(q^{2n} |X|^2 + (1{-}q^{2n})(q^{2n+2} Y^2+c)\big)\zeta_0\nonumber\\
&= v^n \big(q^{2n} (1{-}q^2)(Y^2_0+c)\zeta_0 
+ (1{-}q^{2n})(q^{2n+2} Y^2_0 +c)\zeta_0\big)\nonumber\\
&=\lambda_{n+1}^2 (q^{2n} Y^2_0 +c)\zeta_n
\end{align}
for $\zeta_n\in\cK^s_n$. We assumed above that the pointwise commuting 
self-adjoint operator $|X|^2$ and $Y$ strongly commute. Hence their 
reduced parts on each subspace $\cK^s_n$ strongly commute. 
Therefore, (\ref{x2rel}) implies that 
$|X| \zeta_n=\lambda_{n+1} (q^{2n} Y^2_0 + c)^{1/2}\zeta_n$. 
Recall that $X=v |X|$ and $X^\ast = |X| v^\ast$.  
Renaming $\cK^s_n$ by $\cK_n$ and summarizing the preceding, 
we obtain the following form of the operators $X$, $X^\ast$ and $Y$:
\begin{align*}
&X \zeta_n = \lambda_{n+1} (q^{2n} Y^2_0 + c)^{1/2} \zeta_{n+1},\ \;
X^\ast \zeta_n = \lambda_n (q^{2n-2} Y^2_0 + c)^{1/2} \zeta_{n-1},\\
&Y \zeta_n = q^{2n} Y_0 \zeta_n\quad {\rm on}\ \,
\cK = \oplus^\infty_{n=0} \cK_n,\ \, \cK_n=\cK_0,
\end{align*}
where $Y_0$ is self-adjoint operator with trivial kernel 
on the Hilbert space $\cK_0$. Conversely, it is easy to check 
that these operators $X,X^\ast, Y$ satisfy the relations (\ref{xyrel}), 
so they define indeed a Hilbert space representation 
of the $\ast$-al\-ge\-bra $\cY_c$.
%
\section{Representations of the cross product algebras: \\ 
First approach}
\label{sec-4}
First let us review the representations of the $\ast$-al\-ge\-bra 
$\cO (S^2_{qc})$ from \cite{P}. 
For $c\in [0,\infty)$ and $n\in\dN_0$ 
we abbreviate
$$
\lambda_\pm := 1/2 \pm (c+ 1/4)^{1/2},\ \ 
c_\pm (n) := (c+ \lambda_\pm q^{2n} - (\lambda_\pm q^{2n})^2)^{1/2}.
$$
Let $\Hh^0,\Hh^+_0$ and $\Hh^-_0$ be Hilbert spaces and let $u$ be 
a unitary operator on 
$\Hh^0$. Let $\Hh^\pm = \oplus^\infty_{n=0}\Hh^\pm_n$, 
where $\Hh^\pm_n = \Hh^\pm_0$, and let 
$\Hh = \Hh^0\oplus\Hh^+ \oplus\Hh^-$ for $c\!\in\! (0,+\infty]$ and  
$\Hh = \Hh^0 \oplus \Hh^+$ for $c=0$. The generators of $\cO(S^2_{qc})$ act 
on the Hilbert space $\Hh$ by the following formulas:
\begin{align*}
&c \in [0,+\infty):\quad 
A=0,\  B=c^{1/2} u,\  B^\ast = c^{1/2}\  u^\ast \ \ {\rm on}~ \Hh^0,\\
&A \eta_n = \lambda_\pm q^{2n} \eta_n, \ B \eta_n = c_\pm (n) \eta_{n-1},\  
B^\ast \eta_n = c_\pm (n{+}1) \eta_{n+1} \ \ \mbox{on}~ \Hh^\pm.\\
&\\
&c = +\infty: \quad 
 A = 0, \ B=u, \   B^\ast = u^\ast \ \ \mbox{on}~ \Hh_0,\\
&A \eta_n = \pm q^{2n} \eta_n, \ B \eta_n = (1{-}q^{4n})^{1/2} \eta_{n-1},\  
B^\ast \eta_n = (1{-}q^{4(n+1)})^{1/2} \eta_{n+1}\ \  \mbox{on}~ \Hh^\pm.
\end{align*}
Remind that for $c=0$ there is no Hilbert space $\Hh^-$.
From Proposition 4 in \cite{P} it follows that up to unitary equivalence 
any $\ast$-re\-pre\-sen\-ta\-tion of $\cO(S^2_{qc})$ is of the above form. 
Note that all operators are bounded and $\Hh^0=\ker A$, $A>0$ on $\Hh^+$ 
and $A<0$ on $\Hh^-$.

Next we show that for any $\ast$-rep\-re\-sen\-ta\-tion 
of the algebras $\cO(S^2_{qc}) \rti \cU_q (\mathrm{su}_2)$ or       
$\cO(S^2_{qc})\rti \cU^\prime_q(\mathrm{su}_2)$ the space $\Hh^0$ is $\{0\}$, 
so the operator $A$ is invertible. We carry out the reasoning 
for $\cO(S^2_{qc})\rti \cU_q(\mathrm{su}_2)$. Assuming that the commuting 
self-adjoint operators $A$ and $K$ strongly commute it follows 
that $K$ leaves $\Hh^0 = \ker A$ invariant. 
By (\ref{podrel}) (resp.\ (\ref{prodrel1})), 
$B\Hh^0\subseteq  \Hh^0$. Let $\xi \in \Hh^0$. Using the relation 
$FA = AF {-}q^{-3/2} BK$ we see that 
$$
q^{-3} \| BK \xi\|^2= \|AF \xi\|^2 = \langle q^{-3/2} BK \xi, 
AF\xi\rangle = \langle q^{-3/2} ABK \xi, F\xi\rangle =0.
$$
That is, $BK \xi = 0$ and $AF \xi =0$, so that $F\xi \in\Hh^0$. 
For $c\in (0,+\infty]$, $BK$ is invertible on $\Hh^0$ and hence $\xi=0$. 
For $c=0$ we have $B^\ast \xi=B^\ast F\xi =0$. 
From the cross relation 
$FB^\ast = q^{-1} B^\ast F +q^{-1/2} (1{+}q^2)AK -q^{-1/2} K$ 
we get $K \xi=0$ and so $\xi=0$. Thus, $\Hh^0=\{0\}$.

Consider now a $\ast$-rep\-re\-sen\-ta\-tion of the $\ast$-al\-ge\-bra 
$\hat{\cO}(S^2_{qc})\rti 
\cU^\prime_q(\mathrm{su}_2)$.  Its restriction to $\cO(S^2_{qc})$ is a 
$\ast$-rep\-re\-sen\-ta\-tion of the form described above with  
$\Hh^0 = \{0\}$. 
All operators of the $\ast$-sub\-al\-ge\-bra $\cY_c$ commute with $A$ and $B$. 
Let us assume that the spectral projections of the self-adjoint 
operator $A$ commute also with all operators of $\cY_c$. 
(Note that for a $\ast$-rep\-re\-sen\-ta\-tion by bounded 
operators on a Hilbert space 
the latter fact follows. For unbounded $\ast$-rep\-re\-sen\-ta\-tion 
it does not 
and we restrict ourselves to the class of well-behaved representations 
which satisfy this assumption.) Since $\Hh^\pm_n$ is the eigenspace 
of $A$ at eigenvalue $\lambda_\pm q^{2n}$ and these eigenvalues 
are pairwise distinct,  the operators of $\cY_c$ leave $\Hh^\pm_n$ 
invariant. That is, 
we have a $\ast$-rep\-re\-sen\-ta\-tion of $\cY_c$ on $\Hh^\pm_n$. 
But $\cY_c$ commutes also with $B$. Since $B$ is a weighted shift 
operators with weights $c_\pm (n)\ne 0$ for $n\in\dN$, it follows 
that the representations of $\cY_c$ on $\Hh^\pm_n \!=\! \Hh^\pm_0$ are 
the same for all $n\in\dN_0$. Using the structure of the representation 
of the $\ast$-al\-ge\-bra $\cY_c$ on the Hilbert space $\cK := \Hh^\pm_0$ 
derived in Section 3 and inserting the formulas 
for $X,X^\ast,Y$ and $A,B,B^\ast$ into (\ref{efdef}), 
one obtains the action of the generators 
$e,f,k,k^{-1}$ of $\cU^\prime_q(\mathrm{su}_2)$. 
We do not list these formulas, but we will do so below for the generators 
of $\cU_q(\mathrm{su}_2)$. 

Let us turn to a $\ast$-rep\-re\-sen\-ta\-tion of the larger 
$\ast$-al\-ge\-bra 
$\hat{\cO}(S^2_{qc})\rti \cU_q(\mathrm{su}_2)$. Then, in addition 
to the considerations of the preceding paragraph, we have to deal 
with the generator $K$. We assumed above that $K$ and $A$ are strongly 
commuting self-adjoint operators. Therefore, $K$ commutes with the spectral 
projections of $A$. Hence each space $\Hh^\pm_n$ is reducing for $K$. 
The  relation $BK=q KB$ implies that there is an invertible self-adjoint 
operator $K_0$ on $\Hh^\pm_0$ such that  $K \eta_n = q^n K_0\eta_n$ 
for $\eta_n\in \Hh^\pm_n$. Recall that we have $XK= q KX$ and $YK=KY$ in 
the algebra $\hat{\cO}(S^2_{qc})\rti \cU_q(\mathrm{su}_2)$. 
Inserting for $X$ and $Y$ 
the corresponding operators on 
$\cK\equiv \Hh^\pm_0 = \oplus^\infty_{n=0} \cK_n$ from Section 3 
we conclude that there exists an invertible self-adjoint 
operator $H$ on $\cK_0$ such that $K_0\zeta_n =q^{-n} H\zeta_n$ for 
$\zeta_n\in\cK_n$. Further, since $Y=q K^{-2} A$, 
we have $Y \zeta_0=Y_0\zeta_0=q H^{-2} \lambda_\pm \zeta_0$ 
for $\zeta_0\in\cK_0$.
Inserting the preceding facts and the results from Section 3 
into (\ref{ferel})
and renaming $\cK_0$ by $\cG$
we obtain the following $\ast$-re\-pre\-sen\-tations 
of the cross product 
$\ast$-al\-ge\-bra $\hat{\cO}(S^2_{qc})\rti \cU_q(\mathrm{su}_2)$ 
for $c\in(0,+\infty)$:
\begin{align*}
(I)_{\pm,H}:\quad &A\eta_{nj} =~ \lambda_\pm q^{2n} \eta_{nj}, \ \,
B\eta_{nj} = c_\pm (n) \eta_{n-1,j}, \ \,
B^\ast \eta_{nj} = c_\pm (n{+}1) \eta_{n+1,j},\\
&E \eta_{nj} = ~q^{-1/2} \lambda^{-1} [q^{-n} \lambda_j 
(\lambda_\pm^{-2} q^{-2j} c + H^{-4})^{1/2} H \eta_{n,j-1}\\
&\qquad \quad -q^{-j} (\lambda^{-2}_\pm q^{-2n-2} c 
+ \lambda^{-1}_\pm -q^{2n+2})^{1/2} H \eta_{n+1,j}],\\
&F \eta_{nj} = ~q^{-1/2} \lambda^{-1} [q^{-n} \lambda_{j+1} 
(\lambda^{-2}_\pm q^{-2j-2} c + H^{-4})^{1/2} H \eta_{n,j+1}\\
&\qquad \quad -q^{-j} (\lambda^{-2}_\pm q^{-2n} c 
+ \lambda^{-1}_\pm - q^{2n})^{1/2} H \eta_{n-1,j}],\\
&K \eta_{nj} =~q^{n-j} H \eta_{nj},
\end{align*}
where $H$ is an invertible self-adjoint operator on a Hilbert space $\cG$. 
The representation space is the direct 
sum $\Hh=\oplus^\infty_{n,j=0} \Hh_{nj}$, where $\Hh_{nj} =\cG$.
In the case $c=0$ there is only the representation $(I)_{+,H}$. 
The case $c=+\infty$ has already been treated in [SW].
%
\section{Representations of the cross product algebras: \\ Second approach}
For $l\in\frac{1}{2}\dN_0$ let $T_l$ denote the type 1 spin $l$ 
representations of $\cU_q(\mathrm{su}_2)$. Recall that $T_l$ is 
an irreducible $\ast$-rep\-re\-sen\-ta\-tion of the 
$\ast$-al\-ge\-bra $\cU_q(\mathrm{su}_2)$ 
acting on a $(2l+1)$-dimensional Hilbert space with orthonormal 
basis $\{v^l_j\,;\,j=-l,-l+1,{\dots},l\}$ by the formulas 
(see, for instance, \cite{KS}, p.\ 61)
\begin{equation}\label{kefop}
K v^l_j = q^j v^l_j,\ E v^l_j = [l{-}j]^{1/2} [l{+}j{+}1]^{1/2} v^l_{j+1},\
F v^l_j = [l{-}j{+} 1]^{1/2} [l{+}j]^{1/2} v^l_{j-1}.
\end{equation}
Let $l_0\in \frac{1}{2}\dN_0$ be fixed. The aim of this section 
is classify all $\ast$-rep\-re\-sen\-ta\-tions of the cross product 
$\ast$-al\-ge\-bra  $\cO(S^2_{qc})\rti \cU_q(\mathrm{su}_2)$ having the 
following property:
\medskip

{\it The restriction to $\cU_q(\mathrm{su}_2)$ is the direct 
sum representation $\oplus^\infty_{n=0} T_{l_0+n}$.} 
\hfill $(\ast)$
\medskip\\
In this section 
we shall use 
the generators $x_{j}$ defined by (\ref{xgener}) 
(resp.\ (\ref{xgener1})) rather than $A,B,B^\ast$. 

Suppose we have such a representation. First note that then there exist 
complex numbers $\alpha^\epsilon (l,j)$ and 
$\beta^\epsilon (l,j)$, $\epsilon \!=\! +,0,-$, $l\ge l_0$, 
$j=-l,-l+1,{\dots},l$, $l-l_0\in \dN_0$, such that
\begin{align}\label{xrel1}
x_1 v^l_j &= \alpha^+ (l,j) v^{l+1}_{j+1} + \alpha^0 (l,j) v^l_{j+1} 
+ \alpha^- (l,j) v^{l-1}_{j+1}\\
\label{xrel2}
x_0 v^l_j &= \beta^+ (l,j) v^{l+1}_j + \beta^0 (l,j) v^l_j + 
\overline{\beta^+ (l\!-\!1,j)} v^{l-1}_j\\
\label{xrel3}
x_{-1} v^l_j &= -q^{-1} 
\big(\,\overline{\alpha^- (l\!+\!1, j\!-\!1}) v^{l+1}_{j-1} + 
\overline{\alpha^0(l,j\!-\!1)} v^l_{j-1} +
\overline{\alpha^+(l\!-\!1,j\!-\!1)} v^{l-1}_{j-1}\big),
\end{align}
where the vectors $v^l_{l+1}$, $v^l_{-l-1}$ and the corresponding 
coefficients are set zero. In order to prove this fact 
we argue 
as in \cite{DS}; 
a similar argument appeared also in \cite{S}, p.\ 258. 
Since $Kx_1=q x_1 K$, $x_1v^l_j$ is a weight vector with weight $j+1$. 
Moreover, since 
$E^{l-j+1} x_1 v^l_j=q^{-l+j-1} x_1 E^{l-j+1} v^l_j=0$, $x_1 v^l_j$ is in 
the linear span of vectors  $v^r_{j+1}$, where $r\le l+1$. Similarly, 
replacing $x_1$ by $x_{-1}$ and $E$ by $F$ we conclude that $x_{-1} v^l_j$ 
belongs to the span of vectors  $v^r_{j-1}$, $r\le l+1$. Therefore, 
since $x_{-1}=-q^{-1} x^\ast_1$, $x_{\pm 1} v^l_j\in 
\Lin \{v^{l-1}_{j\pm 1}, v^l_{j\pm 1}, v^{l+1}_{j\pm 1} \}$, 
so $x_{\pm 1} v^l_j$ is of the form (\ref{xrel1}) (resp.\ (\ref{xrel3})). 
From the last two relations of (\ref{podrel}) (resp.\ (\ref{prodrel1})) 
it follows that $x_0v^l_j$ 
is of the form (\ref{xrel2}).

Inserting (\ref{kefop}) and (\ref{xrel1})--(\ref{xrel3}) 
into the equation  $E x_1 v^l_j = q^{-1} x_1 E v^l_j$ we get
\begin{align*}
[l{-}j]^{1/2} [l{+}j{+}3]^{1/2} \alpha^+(l,j)
&=q^{-1} [l{-}j]^{1/2} [l{+}j{+}1]^{1/2} \alpha^+(l,j\!+\!1),\\
[l{-}j{-}1]^{1/2} [l{+}j{+}2]^{1/2} \alpha^0(l,j)
&=q^{-1} [l{-}j]^{1/2} [l{+}j{+}1]^{1/2} \alpha^0(l,j\!+\!1),\\
[l{-}j{-}2]^{1/2} [l{+}j{+}1]^{1/2} \alpha^-(l,j)
&=q^{-1} [l{-}j]^{1/2} [l{+}j{+}1]^{1/2} \alpha^-(l,j\!+\!1).
\end{align*}
The solutions of these recurrence relations are given by 
\begin{align*}
\alpha^+ (l,j) &= q^{-l+j} [l{+}j{+}1]^{1/2} 
[l{+}j{+}2]^{1/2} [2l{+}1]^{-1/2} [2l{+}2]^{-1/2} \alpha^+(l,l),\\
\alpha^0 (l,j) &= q^{-l+j+1} [l{-}j]^{1/2} [l{+}j{+}1]^{1/2} [2l]^{-1/2} 
\alpha^0(l,l\!-\!1),\\
\alpha^- (l,j) &= q^{-l+j+2} [l{-}j{-}1]^{1/2} [l{-}j]^{1/2} [2]^{-1/2} 
\alpha^-(l,l\!-\!2).
\end{align*}
Similarly, from the equation $E x_0 v^l_j=(x_0E+[2]^{1/2} x_1 K) v^l_j$ 
we obtain
\begin{align*}
&[l{-}j{+}1]^{1/2} [l{+}j{+}2]^{1/2} \beta^+(l,j)
\!=\![l{-}j]^{1/2} [l{+}j{+}1]^{1/2} 
\beta^+(l,j{+}1)\!+\![2]^{1/2} q^j \alpha^+(l,j),\\
&[l{-}j]^{1/2} [l{+}j{+}1]^{1/2} \beta^0(l,j)=[l{-}j]^{1/2} [l{+}j{+}1]^{1/2} 
\beta^0(l,j\!+\!1)+[2]^{1/2} q^j \alpha^0(l,j).
\end{align*}
The equation $E x_0 v^l_l =(x_0 E+[2]^{1/2} x_1 K)v^l_l$ yields in addition 
\begin{equation}                                      \label{baplus}
\beta^+ (l,l)= q^l[2]^{1/2} [2l{+}2]^{-1/2} \alpha^+(l,l).
\end{equation}
Further, the equation 
$0=\langle v^l_l, (x_1F+q[2]^{1/2} x_0 K-qF x_1)v^l_l\rangle$ implies 
that 
\begin{equation}\label{banull}
\alpha^0 (l,l\!-\!1)= -[2]^{1/2} [2l]^{-1/2} q^{l+1} \beta^0(l,l).
\end{equation}
Using (\ref{baplus}) and (\ref{banull}) it follows that the above 
recurrence relations for $\beta^+(l,j)$ and $\beta^0(l,j)$ have the 
following solutions:
\begin{align}    \label{bplus}
\beta^+ (l,j) &= 
q^j[l{-}j{+}1]^{1/2} [l{+}j{+}1]^{1/2}[2]^{1/2}[2l{+}1]^{-1/2} [2l{+}2]^{-1/2} 
\alpha^+(l,l),\\
\beta^0 (l,j)&=(1-q^{l{+}j{+}1}[l{-}j][2][2l]^{-1})\beta^0(l,l).\label{bnull}&
\end{align}
From $0=\langle v^l_{l-1}, 
(E x_{-1}-qx_{-1} E-[2]^{1/2} x_0 K)v^{l-1}_{l-1}\rangle$ 
we derive that $\beta^+(l\!-\!1,l\!-\!1)=
-q^{-l}[2l{-}1]^{1/2} \overline{\alpha^-(l,l\!-\!2)}$. 
Combining this equation  with (\ref{baplus}) gives
\begin{equation}\label{aplusminus}
\alpha^- (l,l\!-\!2)= 
-q^{2l-1} [2]^{1/2}[2l{-}1]^{-1/2} [2l]^{-1/2}
\overline{\alpha^+(l\!-\!1,l\!-\!1)}.
\end{equation}
From (\ref{baplus})--(\ref{aplusminus}) we conclude that the representation 
is completely described if the numbers $\alpha^+(l,l)$ and $\beta^0(l,l)$ 
for $l=l_0,l_0+1,{\dots}$ are known. Our next aim is to determine 
these numbers.

Since $\Ker B=\{0\}$ (cf.\ Section \ref{sec-4}), we have 
$\alpha^+(l,l)\neq 0$. By applying a unitary transformation, we can assume 
without loss of generality that $\alpha^+(l,l)$ is a positive real number. 
Furthermore, $\beta^0(l,l)$ is real since $x_0$ is self-adjoint. 
As a consequence, all coefficients occurring 
in (\ref{xrel1})--(\ref{xrel3}) are real. 

We now compute $\beta^0(l_0,l_0)$ for $l_0>0$. Assume first that $c<\infty$. 
We abbreviate $\rho:=1+(q+q^{-1})^2c$. 
The two relations $BB^\ast=q^2 A-q^4 A^2+c$ and $B^\ast B=A-A^\ast +c$ 
lead to the equations 
\begin{align*}
&(1+q^2) (\alpha^0 (l_0,l_0)^2 + \alpha^-(l_0\!+\!1,l_0\!-\!1)^2)= \\
&\qquad\qquad\qquad\qquad\qquad\quad
q^2 \rho + (1-q^2) \beta^0 (l_0, l_0)
-(\beta^0(l_0,l_0)^2 +\beta^+(l_0,l_0)^2),\\
&(1+q^2) \alpha^+ (l_0,l_0)^2 = 
q^2 \rho - (1-q^2) q^2 \beta^0 (l_0, l_0)
-q^4 (\beta^0 (l_0,l_0)^2+\beta^+(l_0,l_0)^2).
\end{align*}
By inserting (\ref{baplus}), (\ref{banull}) and (\ref{aplusminus}) 
into these two relations, we get after some computations
\begin{align}\label{alpharel}
&(1+q^2)q [2l_0{+}3][2l_0{+}2]^{-1} \alpha^+ (l_0,l_0)^2 
=    \nonumber\\
&\qquad\qquad\qquad\qquad\qquad\qquad
q^2 \rho - (1-q^2)q^2\beta^0 (l_0, l_0)-q^4 \beta^0 (l_0,l_0)^2,\\
&(1 + q^2) q^{2l_0+1} [2l_0 {+} 3] [2l_0 {+} 2]^{-1} [2l_0{+}1]^{-1} 
\alpha^+ (l_0,l_0)^2 =\nonumber\\
&\qquad\qquad\quad
q^2 \rho + (1 - q^2) \beta^0 (l_0, l_0) -q^2 [2]^{-1} [2l_0]^{-1} 
\big([2l_0 {+} 3] + q^{2l_0}\big) \beta^0 (l_0,l_0)^2.\nonumber
\end{align}
Eliminating $\alpha^+(l_0,l_0)$ from these equations we obtain 
the following quadratic equation for the coefficient 
$\beta^0(l_0,l_0)$:
$$
q^2\big([2l_0{+}2][2l_0]^{-1} \beta^0(l_0,l_0)\big)^2
-(1-q^2)[2l_0{+}2][2l_0]^{-1} 
\beta^0(l_0,l_0)-\rho=0.
$$
This equation has the solutions $\beta^0(l_0,l_0)_\pm$ given by
\begin{equation}\label{bsol}
\beta^0(l_0,l_0)_\pm=[2l_0][2l_0{+}2]^{-1}(q^{-2} \lambda_\pm -\lambda_\mp).
\end{equation}
In the case $c=\infty$ one gets by analogous calculations the quadratic 
equation $q^2([2l_0{+}2][2l_0]^{-1} \beta^0(l_0,l_0))^2-(q+q^{-1})^2=0$ with 
the solutions 
\begin{equation}                                              \label{bsol1}
\beta^0(l_0,l_0)_\pm=\pm[2l_0][2l_0{+}2]^{-1}(q^{-2}+1).
\end{equation}

Note that $\beta^0(l_0,l_0)_+>0$ and $\beta^0(l_0,l_0)_-<0$ for $l_0>0$. 
Now let $l_0=0$. Then $h(\cdot)=\langle\cdot v^0_0,v^0_0\rangle$ is 
the unique $\cU_q(\mathrm{su}_2)$-invariant state on $\cO(S^2_{qc})$. 
Hence $\beta^0(0,0)= \langle x_0 v^0_0,v^0_0\rangle = h(x_0)=0$. 
That is, formulas (\ref{bsol}) are also valid for $l_0=0$.

Next we determine the coefficient 
$\beta^0(l_0\!+\!n, l_0\!+\!n)$ for $n\in \dN$. 
Using the fact that 
$\|x^n_1 v^{l_0}_{l_0}\|^{-1} x^n_1 v^{l_0}_{l_0}= v^{l_0+n}_{l_0+n}$ 
and the commutation relations of the generators of the algebras 
$\cO(S^2_{qc})$ we derive
\begin{align}\label{x1rel1}
&\| x^n_1 v^{l_0}_{l_0}\|^2 \beta^0 (l_0\!+\!n,l_0\!+\!n)
=\langle x^n_1 v^{l_0}_{l_0}, x_0 x^n_1 v^{l_0}_{l_0}\rangle\nonumber\\
&=\langle x^n_1 v^{l_0}_{l_0},
\big(q^{2n}x^n_1 x_0+(1-q^{2n}) x^n_1\big) v^{l_0}_{l_0}\rangle\nonumber\\
&=\| x^n_1 v^{l_0}_{l_0}\|^2 \big(1-q^{2n}+q^{2n}\beta^0(l_0,l_0)\big)
+q^{2n}\beta^+(l_0,l_0)\langle x^n_1 v^{l_0}_{l_0}, 
x^n_1 v^{l_0+1}_{l_0}\rangle.
\end{align}
On the other hand, applying the relation 
$x_1 v^{l_0}_{l_0}=\alpha^+ (l_0,l_0)v^{l_0+1}_{l_0+1}$ and the 
commutation rules of the cross product algebra we compute
\begin{align*}                                     
&[2l_0 {+} 2]^{1/2} \langle x^n_1 v^{l_0}_{l_0}, 
x^n_1 v^{l_0+1}_{l_0}\rangle = \langle x^n_1 v^{l_0}_{l_0}, 
x^n_1 F v^{l_0+1}_{l_0+1}\rangle\nonumber\\
&=\langle x^n_1 v^{l_0}_{l_0},
\big(q^n F  x^n_1 - q [2n][2]^{-1/2} x_0 x^{n-1}_1 K-
\lambda[n][n{-}1][2]^{-1/2} x^{n-1}_1K\big) v^{l_0+1}_{l_0+1}\rangle\nonumber\\
&=\big(-[2n]q^{l_0+2}\beta^0(l_0\!+\!n,l_0\!+\!n)
-\lambda q^{l_0+1}[n][n{-}1]\big)[2]^{-1/2} \alpha^+(l_0,l_0)^{-1}
\|x^n_1 v^{l_0}_{l_0}\|^2.
\end{align*}
Inserting this identity into (\ref{x1rel1})
and using (\ref{baplus}) with $l$ replaced by $l_0$
we obtain after some further computations
$$
\beta^0 (l_0\!+\!n,l_0\!+\!n)=[2l_0{+}2n{+}2]^{-1}
\big([2l_0{+}2]\beta^0(l_0,l_0)
+(q^{-2}-1)[n][2l_0{+}n{+}1]\big).
$$
By (\ref{bsol}) and (\ref{bsol1}), 
this equation can be written as 
\begin{equation}                                        \label{bnullsol}
\beta^0(l,l)_\pm= [2l{+}2]^{-1}\big([2l_0](q^{-2} \lambda_\pm - \lambda_\mp) 
- (1-q^{-2})[l{-}l_0][l{+}l_0{+}1]\big)
\end{equation}
and 
\begin{equation}                                        \label{bnullsol1}
 \beta^0(l,l)_\pm= [2l{+}2]^{-1}\big(\pm q^{-1}[2][2l_0]
- (1-q^{-2})[l{-}l_0][l{+}l_0{+}1]\big)
\end{equation}
for $c<\infty$ and $c=\infty$, respectively. 

Finally, we determine $\alpha^+(l,l)$. Note 
that (\ref{alpharel}) remains valid if $l_0$ is replaced by 
an arbitrary $l\ge l_0$. For $c<\infty$,   Equation (\ref{alpharel})
is equivalent to 
$$
\alpha^+(l,l)^2= 
[2][2l{+}2][2l{+}3]^{-1}
\big(c+1/4-\{(1+q^{-2})^{-1}(1-\beta^0(l,l))-1/2\}^2\big).
$$
Inserting $\beta^0(l,l)_\pm$ from (\ref{bsol}) into this equation 
we obtain two solutions $\alpha^+(l,l)_\pm$ which can be written as 
\begin{align}
\alpha^+(l,l)_\pm =\,&[2]^{1/2}[2l{+}3]^{-1/2}[2l{+}2]^{-1/2} 
\Big( [2l{+}2]^2 (c+1/4)                                 \nonumber\\
&-\{-\lambda/2 [l{-}l_0{+}1][l{+}l_0{+}1]
\pm [2l_0](c+1/4)^{1/2} \}^2\Big)^{1/2}. 
\label{aplussol}
\end{align}
By similar arguments  one gets in the case $c=\infty$ 
\begin{align}
\alpha^+(l,l)_\pm =\,&[2]^{1/2}[2l{+}3]^{-1/2}[2l{+}2]^{-1/2} 
\Big( [2l{+}2]^2                                \nonumber\\
&-\{-\lambda [2]^{-1}[l{-}l_0{+}1][l{+}l_0{+}1]\pm [2l_0] \}^2\Big)^{1/2}. 
\label{aplussol1}
\end{align}

By Equations (\ref{bnullsol})--(\ref{aplussol1}), all coefficients 
$\alpha^\epsilon (l,j)$ and $\beta^\epsilon (l,j), \epsilon =+,0,-$, 
are known, 
so the corresponding representation of the cross product algebra is 
completely described by formulas (\ref{xrel1})--(\ref{xrel3}). 
Thus we have shown that for $l_0\in\frac{1}{2}\dN$ 
there exits at most two$\ast$-rep\-re\-sen\-ta\-tions of 
$\cO(S^2_{qc})\rti \cU_q(\mathrm{su}_2)$ satisfying condition $(\ast)$. 
For $l_0=0$ there is only one such representation. 
From the corresponding formulas it follows easily that these 
representations are irreducible. 

It still remains to verify that formulas (\ref{xrel1})--(\ref{xrel2}) 
define a representation of $\cO(S^2_{qc})\rti \cU_q(\mathrm{su}_2)$. 
For $l_0=0$ 
this is clear, because then $v^0_0$ is a $\cU_q(\mathrm{su}_2)$-invariant 
vector 
and hence we have the Heisenberg representation. For $l_0\in\frac{1}{2}\dN$ 
this can be done by showing that the defining relations of $\cO(S^2_{qc})$ and 
the cross relations of $\cO(S^2_{qc})\rti \cU_q(\mathrm{su}_2)$ are 
satisfied. We have checked this. We omit the details of these 
lengthy and tedious computations.  
In a forthcoming paper we will construct for any $l_0\in\frac{1}{2}\dN$ 
two inequivalent $\ast$-rep\-re\-sen\-ta\-tions of $\cO(S^2_{qc})\rti 
\cU_q(\mathrm{su}_2)$ 
for which $(\ast)$ holds. By the uniqueness proof given above, 
these representations must be given by 
formulas (\ref{xrel1})--(\ref{xrel2}), so these formulas define 
indeed a $\ast$-rep\-re\-sen\-ta\-tion.

Let us summarize the outcome of the preceding considerations: 
\medskip\\
{\it For any $l_0\in\frac{1}{2}\dN$ there exist 
two $\ast$-rep\-re\-sen\-ta\-tions and for $l_0=0$ there exists one
$\ast$-rep\-re\-sen\-ta\-tion of the cross product 
$\ast$-al\-ge\-bra $\cO(S^2_{qc})\rti \cU_q(\mathrm{su}_2)$ such that the 
restrictions 
to the $\ast$-sub\-al\-ge\-bra} $\cU_q(\mathrm{su}_2)$ {\it is the 
representation} 
$\oplus^\infty_{n=0} T_{l_0+n}$. 
\medskip\\
These representations are irreducible. 
They are labeled by the lower 
indices $+$ and $-$ of the coefficients $\alpha^+(l,l)$ and $\beta^0(l,l)$ 
given by (\ref{aplussol}) and (\ref{bnullsol})
(resp.\ (\ref{aplussol1}) and (\ref{bnullsol1})). 
The underlying Hilbert 
space has an orthonormal basis $\{ v^l_j\,;\, 
j\!=\!-l,\!...,l$, $l\in\dN_0$, $l\ge l_0\}$. 
The generators $E,F,K$ of $\cU_q(\mathrm{su}_2)$ act on this 
basis by 
(\ref{kefop}). The action of the generators $x_1,x_0,x_{-1}$ 
of $\cO(S^2_{qc})$ is given by the following formulas:
\begin{align*}
x_1 v^l_j =\  &q^{-l+j} [l{+}j{+}1]^{1/2} [l{+}j{+}2]^{1/2} 
[2l{+}1]^{-1/2} [2l{+}2]^{-1/2} \alpha^+ (l,l)_\pm v^{l+1}_{j+1}\\
           &-q^{j+2} [l{-}j]^{1/2} [l{+}j{+}1]^{1/2} [2]^{1/2} [2l]^{-1} 
\beta^0 (l,l)_\pm v^l_{j+1}\\
           &-q^{l{+}j{+}1} [l{-}j{-}1]^{1/2} [l{-}j]^{1/2} 
[2l{-}1]^{-1/2} [2l]^{-1/2} \alpha^+ (l\!-\!1,l\!-\!1)_\pm v^{l-1}_{j+1},
\end{align*}
\begin{align*}
x_0 v^l_j =\ &q^j [l{-}j{+}1]^{1/2} [l{+}j{+}1]^{1/2} [2]^{1/2} 
[2l{+}1]^{-1/2} [2l{+}2]^{-1/2} \alpha^+ (l,l)_\pm v^{l+1}_j\\
           &+\big(1-q^{l+j+1} [l{-}j][2][2l]^{-1}\big)
            \beta^0 (l,l)_\pm v^l_j\\
           &+q^j [l{-}j]^{1/2} [l{+}j]^{1/2} [2]^{1/2} 
[2l{-}1]^{-1/2} [2l]^{-1/2} \alpha^+ (l\!-\!1, l\!-\!1)_\pm v^{l-1}_j,
\end{align*}
\begin{align*}
x_{-1} v^l_j =\ &q^{l+j} [l{-}j{+}1]^{1/2} [l{-}j{+}2]^{1/2} 
[2l{+}1]^{-1/2} [2l{+}2]^{-1/2} \alpha^+ (l,l)_\pm v^{l+1}_{j-1}\\
               &+q^j [l{-}j+1]^{1/2} [l{+}j]^{1/2} [2]^{1/2} 
[2l]^{-1} \beta^0 (l,l)_\pm v^l_{j-1}\\
               &-q^{-l+j-1} [l{+}j{-}1]^{1/2} [l{+}j]^{1/2} 
[2l{-}1]^{-1/2} [2l]^{-1/2} \alpha^+ (l\!-\!1,l\!-\!1)_\pm v^{l-1}_{j-1}.
\end{align*}
For $l_0>0$ the two representations are not unitarily equivalent. 
For $l_0=0$ they coincide with the Heisenberg representations 
of $\cO(S^2_{qc})\rti\cU_q(\mathrm{su}_2)$.

\mn


\begin{thebibliography}{999} 

\bibitem[DS]{DS}
Dabrowski, L. and A. Sitarz, 
{\it Dirac operator on the standard Podles quantum sphere.} 
math.QA/0209048.

\bibitem[F]{F}
Fiore, G., {\it On the decoupling of the homogeneous and inhomogeneous parts 
in inhomegeneous quantum groups.} J. Phys. A {\bf 35} (2002), 657--678.


\bibitem[KS]{KS}
Klimyk, K. A. and K. Schm\"udgen, 
{\it Quantum Groups and Their Representations.} Springer, Heidelberg, 1997.


\bibitem[KR]{KR}
Kulish, P. P. and N. Yu. Reshetikhin, 
{\it Quantum linear problem for the sine-Gordon equation 
and higher representations.} 
Zap. Nauchn. Sem. L0MI {\bf 101} (1981), 101--110.

\bibitem[LR]{LR}
Lunts, V. A. and A. L. Rosenberg, 
{\it Differential operators on noncomutative rings.} 
Sel. math. {\bf 3} (1997), 335--359.

\bibitem[P]{P}
Podles, P., {\it Quantum spheres.} Lett. Math. Phys. {\bf 14} (1987), 193--202.

\bibitem[S]{S}
Schm\"udgen, K., {\it Commutator representations of differential calculi 
on the quantum groups $SU_q(2).$} J. Geom. Phys., {\bf 31} (1999), 241--264.

\bibitem[SW]{SW}
Schm\"udgen, K. and E. Wagner, 
{\it Hilbert space representations of cross product algebras.} math.QA/0105185,
to appear in J. Funct. Anal.

\bibitem[SF]{SF}
Sz.-Nagy, B. and C. Foias, 
{\it Analyse harmonique des operateurs de l'espace de Hilbert.} 
Academiai Kiado, Budapest, 1979.



\end{thebibliography}
\end{document}